\documentclass[12pt,reqno]{amsart}

\usepackage{amsmath,amsthm}  
\usepackage{graphicx}  
\usepackage{url}
\usepackage{amsfonts,amssymb} 

\usepackage[margin=1.0in]{geometry}

\usepackage{multicol}

\theoremstyle{theorem}

\newtheorem*{theorem*}{Theorem}
\newtheorem*{corollary*}{Corollary}
 
\theoremstyle{definition}

\newtheorem*{remark}{Remark}

\usepackage{url}
\usepackage{wrapfig}
\usepackage{quoting}
\quotingsetup{vskip=-6pt}

\allowdisplaybreaks

\makeatletter
\@addtoreset{footnote}{page}
\makeatother

\makeatletter
\renewcommand*\env@matrix[1][*\c@MaxMatrixCols c]{%
  \hskip -\arraycolsep
  \let\@ifnextchar\new@ifnextchar
  \array{#1}}
\makeatother

\begin{document}

\title[Nilpotents Leave No Trace]
{Nilpotents  Leave No Trace \\ A Matrix Mystery for Pandemic Times}

\author{ Eric L. Grinberg}
\address{University of Massachusetts Boston}
\email{eric.grinberg@umb.edu}

\begin{abstract}
Reopening a cold case,  inspector  Echelon, high-ranking in the Row Operations Center, is searching for a lost linear map, known to be nilpotent. When a partially decomposed matrix is unearthed, he reconstructs its reduced form, finding it singular. But were its roots nilpotent?
\end{abstract}

\maketitle


%
\let\thefootnote\relax 

\section{Early In the Investigation}
\bigskip

\begin{wrapfigure}[9]{r}{6cm}
\vspace{-17pt}
\begin{center}\includegraphics[width=5.5cm] {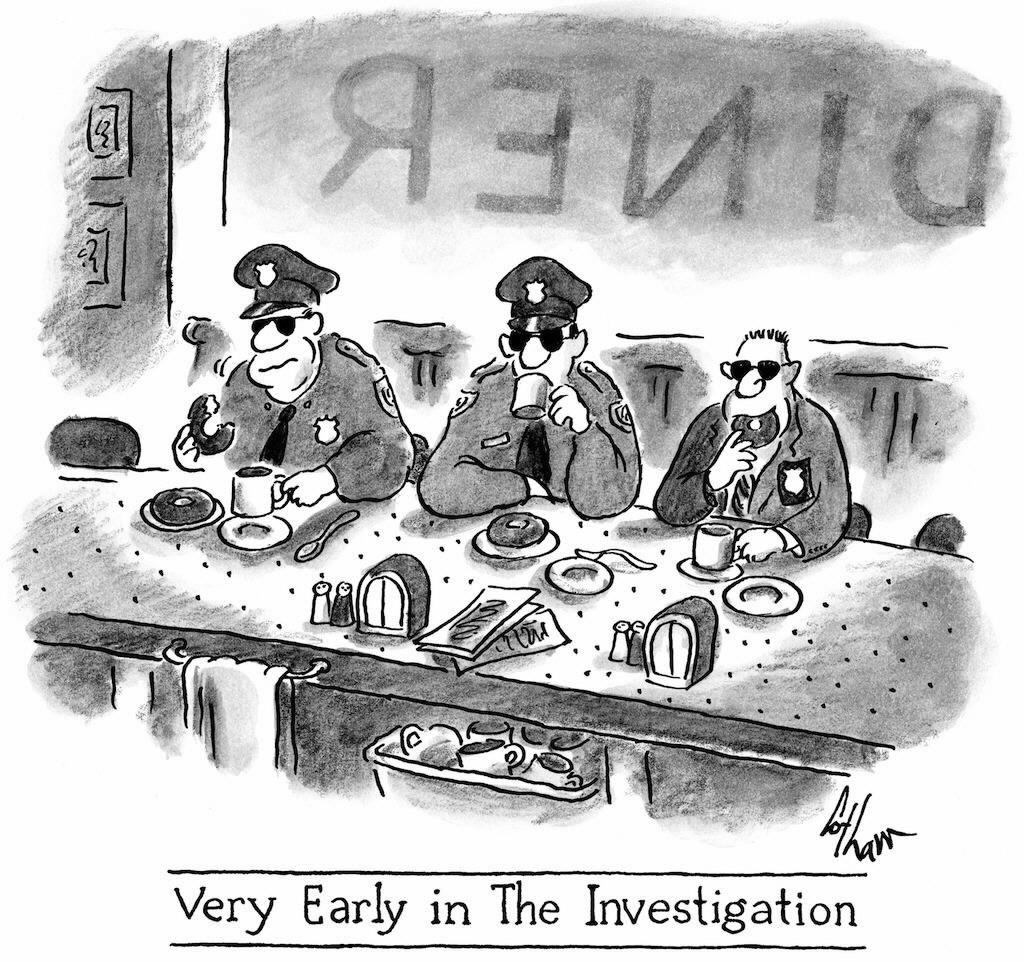} \end{center}
\end{wrapfigure}

In teaching Linear Algebra,  the first topic often is  row reduction \cite{beezer,hoffman-kunze}, including \emph{Row Reduced Echelon Form} (RREF); its applicability is broad and growing. Another topic, surprisingly popular with beginning students, is \emph{nilpotent matrices}.  One naturally wonders about their intersection. For instance, one would expect to find a book exercise asking: 

\smallskip
\qquad \qquad \parbox{7cm}{
\centering
\emph{What can be said about the RREF   \newline
of a nilpotent matrix?}
}

In the early days of the Covid-19 pandemic, as test delivery went remote, demand grew for new, Internet-resistant problems. A limited literature search for the Nilpotent-RREF connection came up short, suggesting potential for take-home final exam questions, hence the note at hand. We'll first explore examples sufficient to settle the $3 \times 3$ case, then consider the general situation. The upshot is the row reduction eliminates all traces of nilpotence.

\section{Stumbling On Evidence}
We refer to  \cite{beezer,beezer-book} for general background on RREF and rank. 
Recall that a square matrix $M$ is \emph{nilpotent} if some power of $M$, say $M^k$, is the zero matrix; the smallest such $k$ is called the \emph{niloptent index} or just \emph{index} of $M$. For instance, the rightmost matrix in \eqref{rref-rank-2} below is nilpotent, of index $3$. Indeed, every \emph{strictly upper-triangular matrix} (square, with zeros on and below the diagonal) is nilpotent. Examining each general type of $3 \times 3$ matrix in RREF, we'll apply row operations to obtain a nilpotent matrix.
\newline

Using the notation of \cite{yutsumura}, the $3 \times 3$ matrices of RREF in rank $1$  are these:

\begin{equation}
\begin{pmatrix} 1 & a & b \\ 0 & 0 & 0 \\ 0 & 0 & 0 \end{pmatrix} , 
\begin{pmatrix} 0 & 1 & c \\ 0 & 0 & 0 \\ 0 & 0 & 0 \end{pmatrix} , \begin{pmatrix} 0 & 0 & 1 \\ 0 & 0 & 0 \\ 0 & 0 & 0 \end{pmatrix}.
\label{rref-rank-1}
\end{equation}

In \eqref{rref-rank-1} the  entries $a,b,c$ are fixed but unspecified and unrestricted constants. 

When working with matrices and their components we'll follow a \emph{left to right} and \emph{up to down} convention. Thus \emph{first} means \emph{leftmost}, etc. We enumerate the rows of a matrix using Roman numberals. Hence \emph{II} is the second row from the top.

The second and third matrices in \eqref{rref-rank-1} are strictly upper triangular, hence nilpotent.
Call  the first matrix $F$. In $F$, if $b=0$, interchange rows $I$ and $III$ to obtain a strictly lower triangular matrix, hence a nilpotent matrix. If $b \ne 0$, perform $III \rightarrow III-\frac{1}{b} I$ (subtract $\frac{1}{b}$ times row $I$ from row $III$ and make that the new row $III$) to obtain
\[
\begin{pmatrix}[rrr] 1 & a & b \\ 0 & 0 & 0 \\
 -\frac{1}{b} & -\frac{a}{b} & -1 \end{pmatrix}.
\]
This matrix squares to zero, hence is nilpotent. 
\newline

Next, we present the RREFs of  $3 \times 3$ matrices with rank $2$:

\begin{equation}
\begin{pmatrix} 1 & 0 & a \\ 0 & 1 & b \\ 0 & 0 & 0 \end{pmatrix},\begin{pmatrix} 1 & a & 0 \\ 0 & 0 & 1 \\ 0 & 0 & 0 \end{pmatrix},\begin{pmatrix} 0 & 1 & 0 \\ 0 & 0 & 1 \\ 0 & 0 & 0 \end{pmatrix}.
\label{rref-rank-2}
\end{equation}

The third matrix is strictly upper triangular, hence nilpotent. For the second matrix, exchange rows $I$ and $III$, then 
perform $I \rightarrow I-(a)II$ to obtain
\[
\begin{pmatrix}[rrr]  0 & 0 & -a \\ 0 & 0 & 1 \\ 1 & a & 0 \end{pmatrix}.
\]
This matrix cubes to zero, hence is nilpotent. As with all $3 \times 3$ matrices of rank $2$, its square does not vanish, so it is nilpotent of index $3$.


So far we have been far from systematic. The leftmost matrix in \eqref{rref-rank-2}, call it $T$, takes a bit more doing. It can be row reduced to the following matrix, which is nilpotent of index $3$:

\begin{equation}
\label{nilpotent-reduction-generic-rank_2}
\begin{pmatrix}[rrr]
-1 & 0 & -a \\
-\frac{b}{a} & 0 & -b \\
-\frac{b - 1}{a} & 1 & 1
\end{pmatrix}.
\end{equation}
We can get from $T$, the leftmost matrix in \eqref{rref-rank-2}
to \eqref{nilpotent-reduction-generic-rank_2} by the following row operations:
\begin{equation}
\label{reduction-steps}
II \leftrightarrow III; \, II \rightarrow II-\frac{b}{a} I; \,  
I \rightarrow (-1) I ; \, III \rightarrow III-\frac{b-1}{b} II.
\end{equation}
This tacitly assumes that $a ,b$ are both nonzero. If both $a$ and $b$ are zero, row swaps turn $T$ into a strictly lower triangular, and hence nilpotent matrix.

If $b=0$  and $a \ne 0$ then 
 \eqref{nilpotent-reduction-generic-rank_2} is still nilpotent and row equivalent to \eqref{nilpotent-reduction-generic-rank_2}, even though the steps we took to get there involve a zero denominator. In case $a=0$, perform row reduction steps 
 \[
 II \rightarrow III ; \, I \rightarrow II; \, II \rightarrow II-bIII,
 \]
 obtaining the following nilpotent matrix:
 \[
 \begin{pmatrix}[rrl]
 0 &   0 & 0 \\
 1 & -b & b^2 \\
 0 &   1 &  b
 \end{pmatrix}.
 \]

 But \eqref{nilpotent-reduction-generic-rank_2} and \eqref{reduction-steps} and all the row manipulations beg the question: how did we come up with these constructs?

\centerline{
}

\section{No Basis For An Investigation}

\begin{quoting}
\flushright 
the facts which you have brought me are so indefinite \\
that we have no basis for an investigation 
\\ 
\smallskip

Sherlock Holmes \\
in \emph{The Adventure of the Dancing Men }
\\

 \bigskip\bigskip

\flushright We have a good working basis, however, on which to start.
\\ \smallskip
Sherlock Holmes \\\
in \emph{A Study In Scarlet }
\end{quoting}

Consider the matrix 
\[
T = 
\begin{pmatrix} 1 & 0 & a \\ 0 & 1 & b \\ 0 & 0 & 0 \end{pmatrix}.
\]
Recall that the (right) \emph{null space} of $T$ is the solution set of the linear system $T \vec v = \vec 0$. One easily checks that the span of the vector $\vec w \equiv \begin{pmatrix} -a & -b & 1 \end{pmatrix}^t$ gives all solutions of this linear system. We will form a basis for $\mathbb R^3$ by extending the one element set $\{ \vec w \}$ and use that to build a nilpotent matrix whose RREF is $T$. Using the familiar notation $\vec e_2 \equiv \begin{pmatrix} 0 & 1 & 0 \end{pmatrix}^t$ and analogs for the canonical basis of $\mathbb R^3$, we write $\vec u \equiv \vec e_2$ and $\vec v = \vec e_3$. Then, if $a$ is a nonzero scalar, $\{\vec u, \vec v, \vec w \}$ is a basis for $\mathbb R^3$. There is a unique  linear transformation $H$ on $\mathbb R^3$ with the properties
$H \vec u = \vec v; \quad H\vec v = \vec w ; \quad H \vec w =0$; we summarize these as follows:
\[
\vec u \longrightarrow \vec v \longrightarrow 
\begin{pmatrix} -a & -b & 1 \end{pmatrix}^t 
\longrightarrow \vec 0.
\]
 (This is not an exact sequence, and not even trying to be one.) Let's find $M$, the matrix representation of $H$. 

The second column of $M$ is the vector $M \vec e_2$, which is already prescribed: it is $\vec e_3$. The third column of $M$ is the vector $M \vec e_3$, prescribed as $\vec w$. What about the first column of $M$? It is $M \vec e_1$, but what's that? We can express $\vec e_1$ in the basis $\{ \vec u , \vec v , \vec w \}$ as follows:
{ \small  \[
(-a) \vec e_1 = 
			 \begin{pmatrix}[r] -a \\ -b \\ 1 \end{pmatrix} 
		-	 \begin{pmatrix}[r]  0 \\  0  \\ 1 \end{pmatrix}
		+	b\begin{pmatrix}[r]  0 \\  1  \\ 0 \end{pmatrix},
\]  }
which we can rewrite as 
\(
\vec e_1 = \frac{-1}{a} \left( \vec w - \vec v +b \vec u \right) .
\)
Thus 
{   \small   \[
M \vec e_1  = 
                 \frac{-1}{a} \left( M\vec w - M\vec v +b M\vec u \right)
              = \frac{-1}{a} \left(  \vec 0  -\vec w +b \vec v \right) 
              = \begin{pmatrix}[r]  -1 \\  -\frac{b}{a}  \\  \frac{1-b}{a} \end{pmatrix},
 \]  }
which matches the first column of \eqref{nilpotent-reduction-generic-rank_2}, and thereby reproduces \eqref{nilpotent-reduction-generic-rank_2}.
\newline
This approach un-begs one question while begging another. The vector space basis procedure here is guaranteed to produce a nilpotent matrix, but how did we know that this matrix will have the requisite RREF, namely $T$? We know that $M$ and $T$ have the same null space: $\text{span} \{ \vec w \}$. We now quote \cite[chp. 2, p. 58] {hoffman-kunze} :

\begin{corollary*}[Hoffman-Kunze]
Let $A$ and $B$  be $m \times n$ matrices over the field $F$. Then $A$ and $B$ are row-equivalent if and only if they have the same row space.
\end{corollary*}
In our context, relating row equivalence to the null space is needed, and such a relation is implicit in the literature, e.g., \cite{hoffman-kunze} again, or \cite[VFSLS]{beezer-book}. A recent posting \cite{gauche-grinberg} gives:

\begin{corollary*}
The null space of a matrix $M$ determines the RREF and the row space of $M$. Hence if two matrices  of the same size have the same null space, they are row equivalent.
\end{corollary*}
 \begin{remark} Everyone knows many famous theorems and some famous lemmas, but there is a dearth of famous corollaries. In fact, the two best known to us are not mathematical, emanating from the Monroe Doctrine.  \end{remark}
 
 \vspace{-10pt}
Thus, since the nilpotent matrix $M$ has the right (right) null space, it has the right RREF as well.

\section{General Impressions}

\begin{quoting}
\flushright 
\hfill Never trust to general impressions, \newline
$\ldots$, but concentrate yourself upon details.
\\ \smallskip
Sherlock Holmes \\
in \emph{A Case of Identity }

I have had no proof yet of the existence
\\ 

\smallskip
Sherlock Holmes \\
in \emph{The Sign Of The Four}
\end{quoting}

\medskip

\begin{theorem*}
Every singular matrix is row equivalent to a nilpotent matrix.
\end{theorem*}
\begin{proof}
Let $M$ be a singular $n \times n$ matrix and take a basis of the (right) null space of $M$, $\{ \vec k_1 , \ldots , \vec k_\ell \}$, where $\ell$ is the nullity of $M$. As $M$ is singular, $\ell$ is greater or equal to $1$. 
If $ \ell =n$ then $M$ is the zero vector, which is nilpotent, and we are done; assume $\ell$ is smaller than $n$. Extend $\{ \vec  k_1, \ldots , \vec k_\ell\}$ to
\( \{ \vec z_1, \ldots \vec z_{n-\ell}, \vec k_1, \ldots , \vec k_\ell \} \),
 a basis of the vector space of all $n \times 1$ columns. We consider a linear map that  annihilates the basis vectors $\vec k_j$ and  ``shifts" each basis vector $\vec z_i$ to the next one, except for $\vec z_{n-\ell}$, which is shifted to $\vec k_1$. This corresponds to a  matrix $N$ with the following properties:
\[
N \vec z_{i-1} =\vec z_i ; \qquad N \vec z_{n-\ell} = \vec k_1,
\]
for all suitable values of $i$. The matrix $N$ is nilpotent of index $n-\ell+1$, with null space spanned by $\{ \vec k_1, \ldots , \vec  k_\ell \}$. Thus $M$ and $N$ share a null space. Hence, by the companion corollary of the \emph{Hoffman-Kunze Corollary}, they have the same RREF and are thereby row equivalent.
\end{proof}


\section{What's It All About? The Aftermath}
Nilpotency figures in the deepest moments of a first course in linear algebra \cite{trotter}. It is particularly accessible to beginning students. Experience indicates that they latch onto to the subject, with  curiosity and enthusiasm; ditto for RREF. Yet, in the literature, the two seldom interact. Why? The theorem may give a clue. In fact, our discussion shows that it's not about nilpotency at all. It's about the null space.
\centerline{
}

\section{Late Inspiration}

In the course of the 2019-2020 academic-pandemic year this author developed the habit of staying up late, delving into the literature, mathematical, fictional and non-fictional. Trying to compose Internet-lookup-resistant take-home final exam questions, he stumbled on the nilpotency-RREF pairing. At the same (late) time, he was reminded of the stories of Jorge Luis Borges, where mathematical points figure into detective stories.  In \emph{Death and The Compass} \cite{borges-aleph,  zalcman88} the mystery hinges on the twists of mathematical reckoning, its inverse and converse, and consequences; an equilateral triangle figures prominently. Recently we noted a new paper about the \emph{Morley triangle} \cite{grinberg-orhon-20}. Morley's is perhaps the most famous of all equilateral triangles. Could a Borges-like story revolve around the Morley triangle, he wondered. In  his half-dormant state he recalled Borges' \emph{The Garden of Forking Paths}, where a secret message is transmitted by the first letter of the name of a person cited prominently in a news story. Could the nilpotency of a matrix serve a similar purpose?  More generally, could a mathematical theorem give rise to a mystery story? (More generally still, is there a functor from the category of mathematics to the category of mystery stories?) That led to the present note. In the vein of the pandemic teachers were asked to exercise particular understanding and accommodation with students. In the same vein one hopes that the reader will do similarly with the would-be\emph{ lockdown literato}  responsible for this pandemic-produced essay.

\underline{Note:} The cartoon included in this essay is by  Frank Cotham. It appeared in the  New Yorker magazine in 2007, and is included in the Cartoon Bank.  Permission for use was obtained by arrangement with  Cond\'e Nast.


\begin{thebibliography}{1}

\bibitem{beezer}
Beezer, R.A.  (2014). Extended Echelon Form and Four Subspaces, \emph{American Math Monthly}, 121:7, 644-647.

\bibitem{beezer-book}
Beezer, R.A.  (2016). \emph{A First Course In Linear Algebra}, Open Source book,
 \url{ http://linear.ups.edu}, accessed May, 2020.
 
\bibitem{borges-lbyrinths}
Borges, J.L.  (1942) . \emph{El jardin de senderos que se bifurcan.} Buenos Aires: Sur. Translated by D. A. Yates (1964), in \emph{Labyrinths: Selected Stories \& Other Writings}, 19-29. New York: New Directions
 
 \bibitem{borges-aleph}
Borges, J.L. (1942). \emph{The Aleph and Other Stories} 1933-1969. NY.: E.P.Dutton, 1970
 
\bibitem{gauche-grinberg} Grinberg, E.L.  (2020) \emph{A Gauche perspective on row reduced echelon form and its uniqueness}, \newline
  \url{arxiv.org/abs/2005.06275}
 
 \bibitem{grinberg-orhon-20} Grinberg E.L., and Orhon, M.  (2020). Morley Trisectors and the Law of Sines with Reflections, 
 \newline \url{arxiv.org/abs/2003.10438}, submitted to \emph{Amer. Math. Monthly}.
 
\bibitem{hoffman-kunze}
Hoffman, K., Kunze, R. (1971). \emph{Linear Algebra}, 2nd ed. Upper Saddle River, NJ:
Prentice-Hall. 

\bibitem{strang-93}
Strang, G. (1993). The Fundamental Theorem of Linear Algebra. \emph{Amer Math Monthly}, 100(9), 848-855. doi:10.2307/2324660.

\bibitem{strang-14}
Strang, G. (2014) The Core Ideas in Our Teaching \emph{Notices  Am. Math.} 
Soc. 61, no. 10, 1243-1245, doi:10.1090/noti1174.

\bibitem{trotter}
Trotter, H. F.  (1961). A Canonical Basis for Nilpotent Transformations. \emph{Amer. Math. Monthly} 68:8, 779-780.

\bibitem{yutsumura}
Tsumura, Y. Problems in Mathematics web site, accessed May 2020,
\newline \url{yutsumura.com/tag/reduced-row-echelon-form/}~ 

\bibitem{zalcman88}
Zalcman, L. (1988). Death and the Calendar, \emph{Hebrew University Studies in
Literature and the Arts} Vol. 16, 97-112.

\end{thebibliography}
\end{document}